\documentclass[11pt,reqno]{amsart}
\usepackage{amscd,amssymb,amsmath,amsthm}
\usepackage{cite}
\newtheorem{thm}[subsection]{Theorem}
\newtheorem{lemma}[subsection]{Lemma}
\newtheorem{pro}[subsection]{Proposition}
\newtheorem{cor}[subsection]{Corollary}

\newtheorem{rk}[subsection]{Remark}
\newtheorem{defn}[subsection]{Definition}

\numberwithin{equation}{section} 

\newcommand{\bea}{\begin{eqnarray}}
\newcommand{\eea}{\end{eqnarray}}

\begin{document}

\title[Rota-Baxter operators]
{Rota-Baxter operators of nilpotent evolution algebras
with maximal nilindex}

\author{Izzat Qaralleh}
\address{Izzat Qaralleh\\
Department of Mathematics\\
Faculty of Science, Tafila Technical
University\\
Tafila, Jordan}
\email{{\tt izzat\_math@yahoo.com}}

\author{Farrukh Mukhamedov}
\address{Farrukh Mukhamedov\\
 Department of Mathematical Sciences\\
College of Science, The United Arab Emirates University\\
P.O. Box, 15551, Al Ain\\
Abu Dhabi, UAE} \email{{\tt far75m@gmail.com} {\tt
farrukh.m@uaeu.ac.ae}}

\author{Otabek Khakimov}
\address{Otabek Khakimov\\
V.I.Romanovskiy Institute  of Mathematics,  Uzbekistan Academy of Sciences\\
Tashkent, Uzbekistan} \email{{\tt khakimovo86@gmail.com}}

\date{Received: xxxxxx; Revised: yyyyyy; Accepted: zzzzzz.
\newline \indent $^{*}$ Corresponding author}

\begin{abstract}
Nilpotent evolution algebras of maximal nilindex admit a natural basis in which the structure matrix is strictly upper triangular. In this paper we classify Rota–Baxter operators of weights zero and one on such algebras. We prove that every Rota–Baxter operator is upper triangular with respect to a natural basis. For weight zero, a strong rigidity phenomenon occurs: the operators are diagonal up to possible perturbations supported in the last basis vector. For weight one, a richer structure appears, including both triangular and non-triangular families, with the diagonal entries governed by a rational recurrence relation. Our results provide a complete description of Rota–Baxter operators on nilpotent evolution algebras of maximal nilindex.\\[2mm]
\end{abstract}


\maketitle


\section{Introduction}

Evolution algebras constitute a class of commutative (generally non-associative) algebras distinguished by the existence of a \emph{natural basis} in which the product of distinct basis vectors vanishes and each square is a linear combination of basis elements.  The algebra is completely encoded by its structure matrix relative to such a basis, and this matrix also defines a directed weighted graph whose vertices correspond to basis elements.  The combinatorial viewpoint makes evolution algebras amenable to methods from graph theory, dynamical systems and genetics \cite{CRR201,CRR202,CNT22}.  In population genetics, for example, evolution algebras were used to model the transmission of types in gametic and zygotic algebras (see Costa's early work on derivations and derivation algebras of gametic and zygotic algebras~\cite{ALP21,Costa1,Costa2,R97,Roz20}).  A comprehensive monograph on evolution algebras and their applications was later given by Tian~\cite{t}, while the survey of Tian and Vojt\u{e}chovsk\'y~\cite{tv} emphasises their role in non-Mendelian genetics.  These algebras sit at the interface between algebraic structures and discrete dynamics; they do not form a variety but their multiplication table satisfies the axioms of a commutative Banach algebra.  As noted by Camacho \emph{et~al.}, evolution algebras are closely related to directed graphs and groups, stochastic processes and mathematical physics~\cite{derevol,QM21,MQ22}.  The interplay between algebraic and dynamical viewpoints has stimulated intensive research on their structural properties and operator theory.

An important structural question concerns nilpotency.  A finite-dimensional evolution algebra $E$ is nilpotent if $E^k=0$ for some $k$; this is equivalent to the absence of oriented cycles in the associated graph, or, equivalently, to the existence of a natural basis in which the structure matrix is strictly upper triangular~\cite{Some_properties}.  Nilpotent evolution algebras have been the subject of several classification efforts.  Two-dimensional evolution algebras were classified over the complex numbers in~\cite{clor} and over the real numbers in~\cite{tv}, and Casado, Molina and Velasco obtained a complete classification of three-dimensional evolution algebras~\cite{3dim}.  In higher dimensions, Elduque and Labra introduced invariants derived from the upper annihilating series to classify nilpotent evolution algebras up to dimension five and constructed families of such algebras using nondegenerate symmetric bilinear forms and commuting diagonalizable endomorphisms~\cite{Elduque}.  Hegazi and Abdelwahab described nilpotent evolution algebras over arbitrary fields and studied their nilpotent index~\cite{heg1}, while Casas, Ladra, Omirov and Rozikov analysed the nilpotent index and dibaricity of evolution algebras and gave conditions for maximal nilindex~\cite{clor,rozomir}.

Among nilpotent evolution algebras, those of \emph{maximal nilindex} play a distinguished role.  An $n$-dimensional evolution algebra has maximal nilindex if $\dim(E^2)=n-1$, or equivalently, if the squares of the natural basis elements span a maximal chain.  In a suitable natural basis, the structure matrix of such an algebra is strictly upper triangular with nonzero superdiagonal entries possibly supplemented by higher-order ``perturbations'' encoded in an index set of off-superdiagonal constants.  Casas \emph{et~al.} showed that nilpotent evolution algebras are nilpotent if and only if their associated graphs contain no oriented cycles, and hence a reordering of a natural basis yields an upper-triangular structure matrix~\cite{Some_properties}.  These extremal nilpotent algebras are analogous to filiform Lie algebras and provide a fertile ground for operator-theoretic investigations. Descriptions of the derivation and automorphisms of evolution algebras with {maximal nilindex} have been investigated in \cite{MKQ24,MKQ20,MKQO19}.

One such operator theoretic constraint is given by Rota--Baxter operators.  For an algebra $A$ over a field $\Bbb K$ and a weight $\theta\in\Bbb K$, a linear map $R\colon A\to A$ is called a Rota--Baxter operator of weight $\theta$ if it satisfies
\begin{equation}
  R(x)R(y)=R\bigl(R(x)y+x\,R(y)+\theta\,xy\bigr),\qquad x,y\in A.
\end{equation}
Rota--Baxter operators generalise integration by parts: when $\theta=0$ and $A$ is the algebra of continuous functions on the real line, the integration operator $P[f](x)=\int_0^x f(t)\,dt$ satisfies $P(f)P(g)=P(P(f)g+f\,P(g))$~\cite{t}.  More generally, Li~Guo's \cite{Guo} exposition on Rota--Baxter algebras emphasises that a Rota--Baxter operator of weight $\lambda$ is a linear map $P$ satisfying $P(x)P(y)=P(xP(y))+P(P(x)y)+\lambda\,P(xy)$ for all $x,y$~\cite{ZSSS}.  Such operators provide a mechanism to construct new algebraic structures from a given one and appear in many areas of mathematics and physics, including probability, combinatorics and renormalization in quantum field theory.  On evolution algebras, derivations and automorphisms have been studied in various contexts~\cite{Ayupov,derevol,Holgate,Gonshor,Costa1,Costa2,GMPQ18,LLR,MQ}, but Rota--Baxter theory on evolution algebras has remained largely unexplored. We notice that Rota--Baxter operators have been recently explored in low dimensional genetic algebras \cite{QMH25,QM24}.

Despite extensive work on derivations, automorphisms, and related linear operators on evolution algebras, Rota–Baxter operators on evolution algebras have received comparatively little attention, especially in higher dimensions and in extremal nilpotent cases. Existing results are mainly confined to low-dimensional or special genetic algebras, and a general classification theory has been lacking.

The goal of this paper is to provide a complete classification of Rota–Baxter operators of weights zero and one on nilpotent evolution algebras of maximal nilindex. Our first main result shows that any Rota–Baxter operator on such an algebra must be upper triangular with respect to a natural basis. This reduces the classification problem to the analysis of polynomial relations among matrix entries, tightly linked to the structure constants of the algebra.

In this paper we undertake the classification of Rota--Baxter operators on nilpotent evolution algebras of maximal nilindex.  Our study reveals a striking dichotomy between weight $\theta=0$ and weight $\theta=1$.  We first show that any Rota--Baxter operator $R=(r_{ij})$ of weight $\theta\in\{0,1\}$ on an $n$-dimensional maximal nilindex evolution algebra is upper-triangular in a natural basis.  This reduces the classification problem to analysing coupled polynomial identities on the coefficients $r_{ij}$.  For weight $\theta=0$ the structure is highly rigid: all off-diagonal entries above the diagonal must vanish except possibly those in the last column, so that $R$ is essentially diagonal plus a single $e_n$-component.  We obtain an explicit description of the space $\mathrm{RB}_0(E)$ of weight-zero Rota--Baxter operators in terms of the index set encoding the off-superdiagonal structure constants.  When this set is empty the operators are diagonal scalings determined by one parameter together with free $e_n$-components up to a threshold index; when long-range structure constants are present, early diagonal entries vanish and the admissible forms become correspondingly restricted.

The situation is more delicate for weight $\theta=1$.  We decompose the space $\mathrm{RB}_1(E)$ into a triangular regime, in which $r_{ij}=0$ for all $i<j<n$, and a non-triangular regime where genuinely upper-triangular interactions occur.  In the triangular regime the diagonal entries of $R$ form an orbit of the rational map $f(x)=\frac{x^2}{2x+1}$, and the index set imposes periodicity constraints determined by the greatest common divisor of certain combinatorial distances.  Over the real numbers, the dynamics of $f$ collapses to explicit normal forms.  In the non-triangular regime the diagonal must be constant (all zeros or all $-1$), and we derive explicit recursive formulae for the remaining coefficients when the index set is empty.  Taken together, our results give a complete classification of Rota--Baxter operators of weights $0$ and $1$ on nilpotent evolution algebras of maximal nilindex.  The contrast between the rigid, diagonal structure in the weight-zero case and the dynamical phenomena emerging for weight one underscores the rich interplay between algebraic, combinatorial and dynamical aspects of evolution algebras.


\section{Evolution algebras}

In this section, we recall  necessary definitions and auxiliary facts about evolution algebras.

Let $\bf{E}$ be a vector space over a field $\mathbb K$.
In what follows, it is assumed that $\mathbb K$ is with the characteristic zero. The vector space $\bf{E}$ with multiplication $\cdot$
is called  an {\it evolution algebra} with respect to a {\it natural basis} $\{{\bf e}_1, {\bf e}_2, . . . \}$ if
the multiplication rule satisfies
$$
{\bf e}_i\cdot {\bf e}_j={\bf 0},\ i\neq j,
$$
$$
{\bf e}_i\cdot{\bf e}_i=\sum_{k}a_{i,k}{\bf e}_k.
$$
From the above definition, it follows that evolution algebras are commutative.
By   $A=(a_{i,j})^n_{i,j=1}$  we denote the matrix of the structural constants
 of the finite-dimensional evolution algebra $\bf{E}$.
Obviously, $\textrm{rank} A =\dim(\bf{E}\cdot\bf{E})$. Hence, for every finite-dimensional evolution algebra the
rank of the matrix does not depend on choice of natural basis.
For an evolution algebra $\bf E$, let us introduce the following sequence
$$
{\bf E}^{k}=\sum_{i=1}^{\lfloor k/2\rfloor}{\bf E}^{i}{\bf E}^{k-i},
$$
where $\lfloor x\rfloor$ denotes the integer part of $x$.

\begin{defn}
An evolution algebra ${\bf E}$ is called nilpotent if there exists some $m\in\mathbb N$
such that ${\bf E}^m=\bf 0$. The smallest $m$ with ${\bf E}^m=\bf 0$ is called the index of nilpotency.
\end{defn}

\begin{thm}\label{thm_ro}\cite{rozomir}
An $n$-dimensional evolution algebra $\bf E$ is nilpotent iff it admits a natural basis such that the matrix of the
structural constants corresponding to $\bf E$ in this basis is represented as follows
$$
\tilde{A}=\left(
\begin{array}{lllll}
0 & \tilde{a}_{1,2} & \tilde{a}_{1,3} & \vdots  & \tilde{a}_{1,n}\\
0 & 0 & \tilde{a}_{2,3} & \vdots  & \tilde{a}_{2,n}\\
\vdots & \vdots & \vdots & \ddots & \vdots\\
0 & 0 & 0 & \vdots & \tilde{a}_{n-1,n}\\
0 & 0 & 0 & \vdots & 0
\end{array}\right).
$$
\end{thm}

\begin{lemma}\label{lemnilEA}
Let $\bf{E}$ be a nilpotent evolution algebra with $\dim({\bf{E}}^2)=n-1$.
Then one can finds a natural basis $\{{\bf e}_1, \dots, {\bf e}_n\}$ such that
\begin{equation}\label{evolalg}
{\bf e}_i^2=\left\{
\begin{array}{lll}
{\bf e}_{i+1}+\sum\limits_{j=i+2}^{n-1}a_{i,j}{\bf e}_j, & i\leq n-3;\\
{\bf e}_{n-1}, & i=n-2;\\
{\bf e}_{n}, & i=n-1;\\
{\bf 0}, & i=n.
\end{array}
\right.
\end{equation}
where $a_{i,j}\in\mathbb K$.
\end{lemma}

\begin{proof} Let ${\bf E}$ be a nilpotent evolution algebra with maximal index of nilpotency, i.e.
$dim({\bf E}^2)=n-1$. Thanks to Theorem \ref{thm_ro} there exists a natural basis
$\{\tilde{\bf e}_1,\dots,\tilde{{\bf e}_n}\}$ such that
$$
\tilde{{\bf e}}_i^2=\left\{
\begin{array}{lll}
\sum\limits_{j=i+1}^{n}\tilde{a}_{i,j}\tilde{{\bf e}}_j, & i\leq n-1;\\
{\bf 0}, & i=n.
\end{array}
\right.
$$
where $a_{i,j}\in\mathbb K$ and $\prod\limits_{i=1}^{n-1}\tilde{a}_{i,i+1}\neq0$.
Without loss of generality we may assume that
$\tilde{a}_{i,i+1}=1$ for any $i<n$. Otherwise instead of the basis $\{\tilde{\bf e}_1,\dots,\tilde{\bf e}_n\}$
we can consider a basis
$$
\tilde{\bf e}_1,\ \tilde{a}_{1,2}\tilde{\bf e}_2,\ \tilde{a}_{1,2}^2\tilde{a}_{2,3}\tilde{\bf e}_3,\ \dots,\ \prod\limits_{k=1}^{n-1}\tilde{a}_{k,k+1}^{2^{n-k-1}}\tilde{\bf e}_n.
$$
Pick the following vectors
$$
{\bf e}_i=
\left\{
\begin{array}{ll}
\tilde{\bf e}_1, & i=1,\\[3mm]
\tilde{\bf e}_i+b_{i}\tilde{\bf e}_n, & 1<i<n,\\[3mm]
\tilde{\bf e}_n, & i=n.
\end{array}
\right.
$$
It is obvious that ${\bf e}_i{\bf e}_j={\bf 0}$ for every $i\neq j$.
Moreover, if we choose coefficients $b_i$ by the following recurrence equation
$$
b_{i+1}=\tilde{a}_{i,n}-\tilde{a}_{i,n-1}b_{n-1}-\dots-\tilde{a}_{i,i+2}b_{i+2},\ \ i<n,
$$
here $b_{n-1}=a_{n-2,n}$, then for ${\bf e}_1,\dots,{\bf e}_n$ one gets
\eqref{evolalg}.
\end{proof}

\begin{rk}\label{remas}
In what follows, we are going to work with nilpotent evolution algebras with the
maximal index of nilpotency.  Therefore, ${\bf E}$ stands for such kind of evolution algebra. Due to Lemma \ref{lemnilEA},
we only consider
evolution algebras whose multiplication table  is given by \eqref{evolalg}.
\end{rk} We identify a set that plays a crucial role in our research.
 For a given matrix $A=(a_{i,j})_{i,j\geq1}^n$  we denote
\begin{equation}\label{A_a_ijneq0}
\mathcal I_A=\{(i,j):i+1<j<n,\ a_{i,j}\neq0\}.
\end{equation}

\section{Rota-Baxter operators}

Let ${\bf E}$ be an algebra over $\mathbb K$ (where $\mathbb K$ is a field of either complex or real numbers) and $\theta\in\mathbb K$.
A linear operator $R$ on the algebra ${\bf E}$
is called {\it a Rota-Baxter operator of weight}
$\theta$ if it satisfies the Rota-Baxter relation of weight
$\theta$:
\begin{equation}\label{RBdef}
R({\bf x})R({\bf y})=R\left(R({\bf x}){\bf y}+{\bf x}R({\bf y})+\theta{\bf x}{\bf y}\right).
\end{equation}
Without loss of generality one can consider the Rota-Baxter operator of weight $\theta\in\{0,1\}$.
By $\mathcal{RB}_\theta({\bf E})$ we denote a set of all Rota-Baxter operators of weight $\theta$ on $\bf E$.

Let ${\bf E}$ be an evolution algebra with natural basis $\{{\bf e}_1,\dots,{\bf e}_n\}$.
Then its Rota-Baxter operator can be represented as follows
\begin{equation}\label{RB1}
R({\bf e}_i)R({\bf e}_i)=R\left(2R({\bf e}_i){\bf e}_i+\theta{\bf e}_i^2\right),\ \ \ i\leq n,
\end{equation}
\begin{equation}\label{RB2}
R({\bf e}_i)R({\bf e}_j)=R\left(R({\bf e}_i){\bf e}_j+{\bf e}_iR({\bf e}_j)\right),\ \ \ i<j\leq n.
\end{equation}
In what follows, we use the matrix representation $R=(r_{i,j})_{i,j=1}^n$, where 
$$
R({\bf e}_i)=\sum\limits_{k=1}^nr_{i,k}{\bf e}_k, \ \ i\in\{1,2,\dots,n\}.
$$

\begin{lemma}\label{lemR-Butri}
Let ${\bf E}$ be an $n$-dimensional evolution algebra with maximal index of nilpotency and
$R=(r_{i,j})_{i,j=1}^n$ be its Rota-Baxter operator of weight $\theta$.
Then $R$ is upper triangular.
\end{lemma}

\begin{proof} Due to $R({\bf e}_n){\bf e}_n={\bf 0}$ we get $R({\bf e}_n)R({\bf e}_n)={\bf 0}$ which is equivalent to
$$
r_{n,1}^2{\bf e}_1^2+r_{n,2}^2{\bf e}_2^2+\dots+r_{n,n-1}^2{\bf e}_{n-1}^2={\bf 0}.
$$
By the linear independence of the vectors ${\bf e}_1^2, {\bf e}_2^2, \dots, {\bf e}_{n-1}^2$
one has
$$
r_{n,1}=r_{n,2}=\dots=r_{n,n-1}=0.
$$
Let us suppose that there exists $i>2$ such that $r_{k,j}=0$ for every $j<k\leq i$. Then by
$R\left(2R({\bf e}_{i-1}){\bf e}_{i-1}+\theta{\bf e}_{i-1}^2\right)=\left(2r_{i-1,i-1}+\theta\right)R({\bf e}_{i-1}^2)$, a multiplication
$R({\bf e}_{i-1})R({\bf e}_{i-1})$ can be written as a linear combination of the vectors
${\bf e}_{i}, {\bf e}_{i+1},\dots, {\bf e}_{n}$.
On the other hand, we have
$$
R({\bf e}_{i-1})R({\bf e}_{i-1})=\sum\limits_{j=1}^{i-2}r_{i-1,j}^2a_{j,j+1}{\bf e}_{j+1}+\sum\limits_{j=i}^n\lambda_j{\bf e}_{j}.
$$
Due to linear independence of ${\bf e}_2,\dots, {\bf e}_{i-1}$ and by
$\prod\limits_{j=1}^{n-1}a_{j,j+1}\neq0$ we infer that
$$
r_{i-1,1}=r_{i-1,2}=\dots=r_{i-1,i-2}=0.
$$
Thus, we have shown upper triangularity of the Rota-Baxter operator.
\end{proof}

\begin{pro}\label{proR-B1}
Let ${\bf E}$ be an $n$-dimensional evolution algebra with maximal index of nilpotency and
$R=(r_{i,j})_{i,j=1}^n$ be its Rota-Baxter operator of weight $\theta$.
Then \eqref{RB1} and \eqref{RB2} can be written as follows

\begin{equation}\label{RB11}
\sum\limits_{k=i}^{m-1}a_{k,m}r_{i,k}^2=(2r_{i,i}+\theta)\sum\limits_{k=i}^{m-1}a_{i,k+1}r_{k+1,m},\ \ \ i<m\leq n.
\end{equation}

\begin{equation}\label{RB21}
\sum\limits_{k=j}^{m-1}r_{i,k}r_{j,k}a_{k,m}=r_{i,j}\sum\limits_{k=j}^{m-1}a_{j,k+1}r_{k+1,m},\ \ \ i<j<m\leq n.
\end{equation}
\end{pro}

\begin{proof}
Let $i<n$. According to Lemma \ref{lemR-Butri}, we get
$R({\bf e}_i)=\sum\limits_{k=i}^nr_{i,k}{\bf e}_k$.
Then keeping in mind ${\bf e}_l{\bf e}_m={\bf 0}$
for any $l\neq m$ we obtain
\begin{equation}\label{eqna1}
R({\bf e}_i)R({\bf e}_i)=\sum\limits_{k=i}^{n-1}r_{i,k}^2{\bf e}_k^2,
\end{equation}
\begin{equation}\label{eqna2}
R(2R({\bf e}_i){\bf e}_i+\theta{\bf e}_i^2)=R\left(2r_{i,i}{\bf e}_i^2+\theta{\bf e}_i^2\right)=(2r_{i,i}+\theta)R({\bf e}_i^2),
\end{equation}
On the other hand, by \eqref{evolalg} we have
\begin{equation}\label{eqna3}
{\bf e}_l^2=\sum\limits_{j=l+1}^na_{l,j}{\bf e}_j,\ \ \ \ \forall l<n.
\end{equation}
Put \eqref{eqna3} into \eqref{eqna1} and \eqref{eqna2} we obtain
\begin{equation}\label{eqna4}
R({\bf e}_i)R({\bf e}_i)=\sum\limits_{k=i}^{n-1}r_{i,k}^2\sum\limits_{m=k+1}^na_{k,m}{\bf e}_m=\sum\limits_{m=i+1}^n\sum\limits_{k=i}^{m-1}a_{k,m}r_{i,k}^2{\bf e}_m
\end{equation}

\begin{equation}\label{eqna5}
R(2R({\bf e}_i){\bf e}_i+\theta{\bf e}_i^2)=(2r_{i,i}+\theta)R({\bf e}_i^2)\sum\limits_{m=i+1}^n\sum\limits_{k=i}^{m-1}a_{i,k+1}r_{k+1,m}{\bf e}_m,
\end{equation}
By equating \eqref{eqna4} and \eqref{eqna5}, we arrive at 
\begin{equation}\label{eqn1}
\sum\limits_{m=i+1}^n\sum\limits_{k=i}^{m-1}a_{k,m}r_{i,k}^2{\bf e}_m=(2r_{i,i}+\theta)\sum\limits_{m=i+1}^n\sum\limits_{k=i}^{m-1}a_{i,k+1}r_{k+1,m}{\bf e}_m,\ \ \ i<n,
\end{equation}
which is equivalent to \eqref{RB11}.
Now we can prove \eqref{RB21}. Since $R({\bf e}_n)=r_{n,n}{\bf e}_n$ and ${\bf e}_i{\bf e}_n={\bf 0}$, $i=\overline{1,n-1}$
one can check that \eqref{RB2} reduces to ${\bf 0}={\bf 0}$ for $j=n$.
Therefore, it is enough to consider $j<n$. Let us take $i\in\{1,\dots,n-1\}$ and $j\in\{i+1,\dots,n-1\}$.
Then 
\begin{equation}\label{eqna6}
R({\bf e}_i)R({\bf e}_j)=\sum\limits_{k=j}^{n-1}r_{i,k}r_{j,k}{\bf e}_k^2,
\end{equation}
\begin{equation}\label{eqna7}
R\left(R({\bf e}_i){\bf e}_j+{\bf e}_iR({\bf e}_j)\right)=R(R({\bf e}_i){\bf e}_j)=r_{i,j}R({\bf e}_j^2).
\end{equation}
Plugging \eqref{eqna3} into \eqref{eqna6} and \eqref{eqna7} we immediately find
\begin{equation}\label{eqna8}
R({\bf e}_i)R({\bf e}_j)=\sum\limits_{m=j+1}^n\sum\limits_{k=j}^{m-1}r_{i,k}r_{j,k}a_{k,m}{\bf e}_m,
\end{equation}
\begin{equation}\label{eqna9}
R\left(R({\bf e}_i){\bf e}_j+{\bf e}_iR({\bf e}_j)\right)=r_{i,j}\sum\limits_{m=j+1}^n\sum\limits_{k=j}^{m-1}a_{j,k+1}r_{k+1,m}{\bf e}_m.
\end{equation}
By equating \eqref{eqna8} and \eqref{eqna9}, one gets
\begin{equation}\label{eqn2}
\sum\limits_{m=j+1}^n\sum\limits_{k=j}^{m-1}r_{i,k}r_{j,k}a_{k,m}{\bf e}_m=r_{i,j}\sum\limits_{m=j+1}^n\sum\limits_{k=j}^{m-1}a_{j,k+1}r_{k+1,m}{\bf e}_m,\ \ \ i<j<n.
\end{equation}
It is clear that \eqref{eqn2} implies \eqref{RB21}.
\end{proof}
In the next sections we are going to examine the cases where $\theta=0$ and $\theta=1$ individually.

\section{Description of $\mathcal{RB}_0({\bf E})$}
Assume that $\theta=0$. In this case equalities \eqref{RB11} and \eqref{RB21} have the following forms

\begin{equation}\label{RB110}
\sum\limits_{k=i}^{m-1}a_{k,m}r_{i,k}^2=2r_{i,i}\sum\limits_{k=i}^{m-1}a_{i,k+1}r_{k+1,m},\ \ \ i<m\leq n.
\end{equation}

\begin{equation}\label{RB210}
\sum\limits_{k=j}^{m-1}r_{i,k}r_{j,k}a_{k,m}=r_{i,j}\sum\limits_{k=j}^{m-1}a_{j,k+1}r_{k+1,m},\ \ \ i<j<m\leq n.
\end{equation}

\begin{lemma}\label{lemR-B2}
Let ${\bf E}$ be an $n$-dimensional evolution algebra with maximal index of nilpotency and
$R=(r_{i,j})_{i,j=1}^n$ be its Rota-Baxter operator of weight $\theta=0$.
Then
\begin{equation}\label{r_diag}
r_{i,i}^2=2r_{i,i}r_{i+1,i+1},\ \ \ \ \forall i<n.
\end{equation}
\end{lemma}

\begin{proof}
Let $i<n$. Then for $m=i+1$ equality \eqref{RB110} can be written as follows
$$
a_{i,i+1}r_{i,i}^2=2r_{i,i}a_{i,i+1}r_{i+1,i+1}.
$$
Then by $a_{i,i+1}=1$ we get \eqref{r_diag}.
\end{proof}

\begin{cor}\label{corR-B3}
Let ${\bf E}$ be an $n$-dimensional evolution algebra with maximal index of nilpotency and
$R=(r_{i,j})_{i,j=1}^n$ be its Rota-Baxter operator of weight $\theta=0$.
Then the followings are true:
\begin{itemize}
\item[$(a)$] if $r_{i,i}=0$ for some $i\in\{2,\dots,n\}$ then $r_{j,j}=0$ for any $j<i$.
\item[$(b)$] if $r_{i,i}\neq0$ for some $i\in\{1,2,\dots,n-1\}$ then $r_{j,j}=2^{i-j}r_{i,i}$ for any $j>i$.
\end{itemize}
\end{cor}
The proof is straightforward since it follows directly from the equalities \eqref{r_diag}.

\begin{lemma}\label{lemR-Bvform}
Let ${\bf E}$ be an $n$-dimensional evolution algebra with maximal index of nilpotency and
$R=(r_{i,j})_{i,j=1}^n$ be its Rota-Baxter operator of weight $\theta=0$.
Then $r_{i,j}=0$ for every $i<j<n$.
\end{lemma}

\begin{proof} By Lemma \ref{lemR-Butri}, Rota-Baxter operator $R$ has upper triangular form, i.e.
$$
R({\bf e}_i)=\sum\limits_{j=i}^nr_{i,j}{\bf e}_j,\ \ \ \ \forall i\leq n.
$$
First, we show that $r_{i,i+1}=0$ for every $i<n-1$. We then conclude the proof by demonstrating that
$\mathcal I_R=\emptyset$.
Let us take any $i\in\{1,2,\dots,n-2\}$.
After rewriting \eqref{RB210} for $j=i+1$ and $m=i+2$ we have
$$
r_{i,i+1}r_{i+1,i+1}a_{i+1,i+2}=r_{i,i+1}a_{i+1,i+2}r_{i+2,i+2},
$$
which together $a_{i+1,i+2}\neq0$ implies
\begin{equation}\label{p1}
r_{i,i+1}r_{i+1,i+1}=r_{i,i+1}r_{i+2,i+2}.
\end{equation}
By \eqref{p1}, $r_{i,i+1}=0$ if $r_{i+1,i+1}\neq r_{i+2,i+2}$.
Now we show that $r_{i+1,i+1}=r_{i+2,i+2}$ implies $r_{i,i+1}=0$.
Assume that
$r_{i+1,i+1}=r_{i+2,i+2}$.
By Lemma \ref{lemR-B2} one has
$$
r_{i+1,i+1}^2=2r_{i+1,i+1}r_{i+2,i+2}.
$$
Hence, $r_{i+2,i+2}=0$. Due to Corollary \ref{corR-B3}, we obtain $r_{k,k}=0$ for every $k\leq i+2$.
Then rewrite \eqref{RB110} for $m=i+2$ we obtain
$$
a_{i+1,i+2}r_{i,i+1}=0.
$$
Since $a_{i+1,i+2}=1$ we infer that $r_{i,i+1}=0$.
Thus we have shown that
$$
R=\left(
\begin{array}{llllllll}
r_{1,1} & 0 & r_{1,3} & r_{1,4} & \ldots & r_{1,n-2} & r_{1,n-1} & r_{1,n}\\
0 & r_{2,2} & 0 & r_{2,4} & \ldots & r_{2,n-2} & r_{2,n-1} & r_{2,n}\\
0 & 0 & r_{3,3} & 0 & \ldots & r_{3,n-2} & r_{3,n-1} & r_{3,n}\\
\vdots & \vdots & \vdots & \vdots & \ddots & \vdots & \vdots & \vdots\\
0 & 0 & 0 & 0 & \ldots & r_{n-2,n-2} & 0 & r_{n-2,n}\\
0 & 0 & 0 & 0 & \ldots & 0 & r_{n-1,n-1} & r_{n-1,n}\\
0 & 0 & 0 & 0 & \ldots & 0 & 0 & r_{n,n}
\end{array}
\right).
$$
Now we check $\mathcal I_R=\emptyset$. Suppose contrary, i.e. $\mathcal I_R\neq\emptyset$.
Then one can define
$$
i_0:=\max\{i: \exists j>i+1\ \mbox{such that}\ (i,j)\in\mathcal I_R\}.
$$
It is obvious that
$$
R({\bf e}_k)=r_{k,k}{\bf e}_k+r_{k,n}{\bf e}_n,\ \ \ \ \forall k\in\{i_0+1,\dots,n-1\}.
$$
Let us pick a $(i_0,j_0)\in\mathcal I_R$.
Rewrite \eqref{RB210} for $i_0$, $j_0$ and $m=j_0+1$
$$
r_{i_0,j_0}r_{j_0,j_0}a_{j_0,j_0+1}=r_{i_0,j_0}a_{j_0,j_0+1}r_{j_0+1,j_0+1}.
$$
Keeping in mind $a_{j_0,j_0+1}=1$ from the last equality
we get
\begin{equation}\label{p2}
r_{i_0,j_0}r_{j_0,j_0}=r_{i_0,j_0}r_{j_0+1,j_0+1}.
\end{equation}
By $r_{i_0,j_0}\neq0$, equality \eqref{p2} yields
$$
r_{j_0,j_0}=r_{j_0+1,j_0+1}.
$$
One the other hand by
Lemma \ref{lemR-B2}, one has $r_{j_0,j_0}^2=2r_{j_0,j_0}r_{j_0+1,j_0+1}$. So, we infer that
$r_{j_0+1,j_0+1}=0$. Thanks to Corollary \ref{corR-B3} one finds $r_{i_0,i_0}=0$. Then rewrite \eqref{RB110}
for $i_0$
\begin{equation}\label{p3}
\left\{
\begin{array}{ll}
a_{i_0,i_0+1}r_{i_0,i_0}^2=0,\\
a_{i_0,i_0+1}r_{i_0,i_0}^2+a_{i_0+1,i_0+2}r_{i_0,i_0+1}^2=0,\\
\dots\dots\dots\dots\dots\dots\dots\dots\dots\dots\\
a_{i_0,i_0+1}r_{i_0,i_0}^2+a_{i_0+1,i_0+2}r_{i_0,i_0+1}^2+\dots+a_{n-1,n}r_{i_0,n-1}^2=0.
\end{array}
\right.
\end{equation}
By $a_{i,i+1}=1$ for every $i<n$, from \eqref{p3} we obtain
$r_{i_0,j}=0$ for any $j>i_0$. In particularly, $r_{i_0,j_0}=0$.
This contradicts to $(i_0,j_0)\in\mathcal I_{R}$.
So, we infer that $\mathcal I_R=\emptyset$.
\end{proof}

\begin{rk}\label{remRB1}
Let ${\bf E}$ be an evolution algebra with structural matrix
$$
A=\left(
\begin{array}{llllllll}
0 & 1 & a_{1,3} & a_{1,4} & \dots & a_{1,n-1} & 0\\
0 & 0 & 1 & a_{2,4} &\dots & a_{2,n-1} & 0\\
0 & 0 & 0 & 1 & \dots & a_{3,n-1} & 0\\
\vdots & \vdots & \vdots & \vdots & \ddots & \vdots & \vdots\\
0 & 0 & 0 & 0 & \dots & 1 & 0\\
0 & 0 & 0 & 0 & \dots & 0 & 1\\
0 & 0 & 0 & 0 & \dots & 0 & 0
\end{array}\right)
$$ in a
natural basis
$\{{\bf e}_i\}_{i=1}^n$ and $R$ be a linear operator on $\bf E$.
If $R$ is Rota-Baxter operator of weight $\theta=0$
then by Lemma \ref{lemR-Bvform}, $R$ can be written as follows
\begin{equation}\label{RBFform}
R({\bf e}_i)=\left\{
\begin{array}{ll}
r_{i,i}{\bf e}_i+r_{i,n}{\bf e}_n, & i<n,\\
r_{n,n}{\bf e}_n, & i=n.
\end{array}
\right.
\end{equation}
It is easy to check that every linear mapping on ${\bf E}$ defined as \eqref{RBFform} satisfies
condition \eqref{RB210}. This means that linear mapping \eqref{RBFform} is Rota-Baxter operator on ${\bf E}$
if and only if
\begin{equation}\label{RB1101}
\left\{
\begin{array}{ll}
r_{i,i}^2=2r_{i,i}r_{i+1,i+1}, & \mbox{if}\ i<n,\\
r_{i,i}^2=2r_{i,i}r_{j,j}, & \mbox{if}\ (i,j)\in\mathcal I_A,\\
r_{i,i}\left(r_{i+1,n}+\sum\limits_{(i,j)\in\mathcal I_A}a_{i,j}r_{j,n}\right)=0, & \mbox{if}\ i<n-2,\\
r_{n-2,n-2}r_{n-1,n}=0.
\end{array}
\right.
\end{equation}
\end{rk}

\begin{thm}\label{thmRB0}
Let ${\bf E}$ be an evolution algebra with structural matrix
$$
A=\left(
\begin{array}{llllllll}
0 & 1 & a_{1,3} & a_{1,4} & \dots & a_{1,n-1} & 0\\
0 & 0 & 1 & a_{2,4} &\dots & a_{2,n-1} & 0\\
0 & 0 & 0 & 1 & \dots & a_{3,n-1} & 0\\
\vdots & \vdots & \vdots & \vdots & \ddots & \vdots & \vdots\\
0 & 0 & 0 & 0 & \dots & 1 & 0\\
0 & 0 & 0 & 0 & \dots & 0 & 1\\
0 & 0 & 0 & 0 & \dots & 0 & 0
\end{array}\right)
$$ in a
natural basis
$\{{\bf e}_i\}_{i=1}^n$ and $R$ be a linear operator on $\bf E$. Then  $R\in\mathcal{RB}_0({\bf E})$ if and only if
it has one of the following forms:
\begin{enumerate}
\item[$(i)$]
if $\mathcal I_A=\emptyset$, then
\begin{enumerate}
\item[$(R_1)$.] there exist $\alpha\in\mathbb K\setminus\{0\}$, $\beta_1\in\mathbb K$ such that
$$
R({\bf e}_i)=
\left\{
\begin{array}{ll}
2^{n-1}\alpha{\bf e}_1+\beta_1{\bf e}_n, & i=1,\\
2^{n-i}\alpha{\bf e}_i, & i>1.
\end{array}
\right.
$$
\item[$(R_k)$.] $(1<k<n)$ there exist $\alpha\in\mathbb K\setminus\{0\}$, $\beta_1,\beta_2,\dots,\beta_{k}\in\mathbb K$ such that
$$
R({\bf e}_i)=
\left\{
\begin{array}{ll}
\beta_i{\bf e}_n, & i<k,\\
2^{n-k}\alpha{\bf e}_k+\beta_k{\bf e}_n, & i=k,\\
2^{n-i}\alpha{\bf e}_i, & i>k.
\end{array}
\right.
$$
\item[$(R_n)$.] there exist $\beta_1,\beta_2,\dots,\beta_{n}\in\mathbb K$ such that
$$
R({\bf e}_i)=\beta_i{\bf e}_n,\ \ \ \ \forall i\leq n.
$$
\end{enumerate}

\item[$(ii)$] if $\mathcal I_A\neq\emptyset$ and $i_0=\max\{i: \exists j>i+1\ \mbox{such that}\ (i,j)\in\mathcal I_A\}$,
then
\begin{enumerate}
\item[$(R_k)$.] $(i_0<k<n)$ there exist $\alpha\in\mathbb K\setminus\{0\}$ and $\beta_1,\dots,\beta_{i_0+1}\in\mathbb K$ such that
$$
R({\bf e}_i)=\left\{
\begin{array}{lll}
\beta_i{\bf e}_n, & i\leq k,\\
2^{n-k}\alpha{\bf e}_{k}+\beta_{k}{\bf e}_n, & i=k,\\
2^{n-i}\alpha{\bf e}_i, & i>k.
\end{array}
\right.
$$
\item[$(R_{n})$.] there exist $\beta_1,\dots,\beta_{n}\in\mathbb R$ such that
$$
R({\bf e}_i)=\beta_i{\bf e}_n,\ \ \ \ \forall i\leq n.
$$
\end{enumerate}
\end{enumerate}
\end{thm}

\begin{proof} Due to Remark \ref{remRB1} every $R\in\mathcal{RB}_0({\bf E})$
has the following form
$$
R({\bf e}_i)=
\left\{
\begin{array}{ll}
r_{i,i}{\bf e}_i+r_{i,n}{\bf e}_n, & i<n,\\
r_{n,n}{\bf e}_n, & i=n.
\end{array}
\right.
$$
Moreover, the defined linear mapping $R$ is a Rota-Baxter operator of ${\bf E}$ for $\theta=0$
if and only if it satisfies
\eqref{RB1101}.

$(i)$ Suppose that $\mathcal I_A=\emptyset$, then it yields that \eqref{RB1101} can be written as follows
\begin{equation}\label{teng1}
\left\{\begin{array}{ll}
r_{i,i}^2=2r_{i,i}r_{i+1,i+1}, & i<n,\\
r_{i,i}r_{i+1,n}=0, & i<n-1.
\end{array}
\right.
\end{equation}
Assume that $r_{1,1}\neq0$. Then from \eqref{teng1} we find
$r_{i,i}=2^{1-i}r_{1,1}$, $i=\overline{2,n}$ and $r_{j,n}=0$, $j=\overline{2,n-1}$.
In this case Rota-Baxter operator has the form $(R_1)$.
Assume that $r_{k,k}\neq0$ and $r_{k-1,k-1}=0$ for some $k\in\{2,\dots,n-1\}$. By
Corollary \ref{corR-B3} we have $r_{i,i}=0$ for every $i<k$ and $r_{j,j}\neq0$ for any
$j\geq k$. In this case the Rota-Baxter operator has the form $(R_k)$.
Finally, consider a case $r_{n-1,n-1}=0$. By Corollary \ref{corR-B3} we arrive at 
$$
r_{1,1}=\dots=r_{n-1,n-1}=0.
$$
One can check that \eqref{teng1} holds for every values
of $r_{i,n}$, $i=\overline{1,n}$. Thus, the Rota-Baxter operator has the form $(R_n)$.

$(ii)$ Let $\mathcal I_A\neq\emptyset$.
We need to check \eqref{RB1101}.
Keeping in mind $\mathcal I_A\neq\emptyset$ one can define
$$
i_0:=\max\{i: \exists j>i+1\ \mbox{such that}\ (i,j)\in\mathcal I_A\}.
$$
It is obvious that $i_0<n-2$. We show that $r_{i,i}=0$ for any $i\leq i_0$.
Let us pick an arbitrary $j_0>i_0+1$ such that $(i_0,j_0)\in\mathcal I_A$.
Then from the second equality of \eqref{RB1101} we get
\begin{equation}\label{eq3}
r_{i_0,i_0}^2=2r_{i_0,i_0}r_{j_0,j_0}.
\end{equation}
Suppose that $r_{i_0,i_0}\neq0$. Then \eqref{eq3} implies $r_{i_0,i_0}=2r_{j_0,j_0}$. 

On the other hand, the first equality of \eqref{RB1101}
yileds that $r_{i_0,i_0}=2^{j_0-i_0}r_{j_0,j_0}$.
Hence, $2^{j_0-i_0}=2$, which contradicts to $j_0>i_0+1$. This means that it holds $r_{i_0,i_0}=0$.
Consequently, Corollary \ref{corR-B3} yields that $r_{i,i}=0$ for every $i\leq i_0$.
Hence,
\begin{equation}\label{teng4}
\left\{
\begin{array}{ll}
r_{i,i}=0, & i\leq i_0,\\
r_{i,i}^2=2r_{i,i}r_{i+1,i+1}, & i_0<i<n,\\
r_{i,i}r_{i+1,n}=0, & i_0<i<n-1.
\end{array}
\right.
\end{equation}
Let $r_{k,k}\neq0$ and $r_{k-1,k-1}=0$ for some $k\in\{i_0+1,\dots,n-1\}$.
Then there exist $\alpha\in\mathbb K\setminus\{0\}$ and $\beta_1,\dots,\beta_{k}\in\mathbb K$
such that
$$
\left\{
\begin{array}{ll}
r_{i,i}=0, & i<k,\\
r_{i,i}^2=2^{n-i}\alpha, & i\geq k,\\
r_{i,n}=\beta_i, & i\leq k,\\
r_{i,n}=0, & k+1<i<n.
\end{array}
\right.
$$
In the considered case, the Rota-Baxter operator has the form $(R_k)$.

Let $r_{n-1,n-1}=0$. Then there exist $\beta_1,\dots,\beta_n\in\mathbb K$ such that
$$
\left\{
\begin{array}{ll}
r_{i,i}=0, & i< n,\\
r_{i,n}=\beta_i, & i\leq n.
\end{array}
\right.
$$
which implies the Rota-Baxter operator has the form $(R_{n})$. This completes the proof.
\end{proof}

\begin{rk} We note that if $n\leq3$ then $\mathcal I_A=\emptyset$. Therefore, by Theorem \ref{thmRB0},
one has
$$
\mathcal{RB}_0({\bf E})=\left\{
\left(
\begin{array}{ll}
2\alpha & \beta\\
0 & \alpha
\end{array}
\right):\ \forall\alpha,\beta\in\mathbb K, \alpha\neq0
\right\}\bigcup
\left\{
\left(
\begin{array}{ll}
0 & \beta\\
0 & \alpha
\end{array}
\right):\ \forall\alpha,\beta\in\mathbb K\right\}
$$
for $n=2$, and
\begin{eqnarray*}
\mathcal{RB}_0({\bf E})
&=&\left\{
\left(
\begin{array}{lll}
4\alpha & 0 & \beta\\
0 & 2\alpha & \gamma\\
0 & 0 & \alpha
\end{array}
\right):\ \forall\alpha,\beta,\gamma\in\mathbb K, \alpha\neq0
\right\}\\
&&\bigcup
\left\{
\left(
\begin{array}{lll}
0 & 0 & \beta\\
0 & 2\alpha & \gamma\\
0 & 0 & \alpha
\end{array}
\right):\ \forall\alpha,\beta,\gamma\in\mathbb K,\ \alpha\neq0\right\}\\
&&\bigcup
\left\{
\left(
\begin{array}{lll}
0 & 0 & \beta\\
0 & 0 & \gamma\\
0 & 0 & \alpha
\end{array}
\right):\ \forall\alpha,\beta,\gamma\in\mathbb K\right\}
\end{eqnarray*}
for $n=3$, respectively.
\end{rk}

\section{On description of $\mathcal{RB}_1({\bf E})$}

Let $\theta=1$. In this setting, the equalities \eqref{RB11} and \eqref{RB21} are reduced to 

\begin{equation}\label{RB111}
\sum\limits_{k=i}^{m-1}a_{k,m}r_{i,k}^2=(2r_{i,i}+1)\sum\limits_{k=i}^{m-1}a_{i,k+1}r_{k+1,m},\ \ \ i<m\leq n.
\end{equation}

\begin{equation}\label{RB211}
\sum\limits_{k=j}^{m-1}r_{i,k}r_{j,k}a_{k,m}=r_{i,j}\sum\limits_{k=j}^{m-1}a_{j,k+1}r_{k+1,m},\ \ \ i<j<m\leq n.
\end{equation}
Thanks to Lemma \ref{lemR-Butri} every $R\in\mathcal{RB}_1({\bf E})$ is upper triangular, i.e.
\begin{equation}\label{Rup}
R=\left(
\begin{array}{llllll}
r_{1,1} & r_{1,2} & r_{1,3} & \ldots & r_{1,n-1} & r_{1,n}\\
0 & r_{2,2} & r_{2,3} & \ldots & r_{2,n-1} & r_{2,n}\\
0 & 0 & r_{3,3} & \ldots & r_{3,n-1} & r_{3,n}\\
\vdots & \vdots & \vdots & \ddots & \vdots & \vdots\\
0 & 0 & 0 & \ldots & r_{n-1,n-1} & r_{n-1,n}\\
0 & 0 & 0 & \ldots & 0 & r_{n,n}
\end{array}
\right).
\end{equation}
Let us denote
$$
\mathcal{RB}_1^\vartriangle({\bf E})=\left\{R\in\mathcal{RB}_1({\bf E}): r_{i,j}=0,\ \mbox{for any}\ i<j<n\right\},
$$
$$
\mathcal{RB}_1^\blacktriangle({\bf E})=\mathcal{RB}_1({\bf E})\setminus\mathcal{RB}_1^\vartriangle({\bf E}).
$$
We will describe the sets $\mathcal{RB}_1^\vartriangle({\bf E})$ and $\mathcal{RB}_1^\blacktriangle({\bf E})$.

\subsection{Case $\mathcal{RB}_1^\vartriangle({\bf E})$}

Note that in this case operator $R$ defined as \eqref{Rup} satisfies \eqref{RB211}. The condition \eqref{RB111}
is equivalent to the following
\begin{equation}\label{RB1110}
\left\{
\begin{array}{ll}
r_{i,i}^2=(2r_{i,i}+1)r_{i+1,i+1}, & \mbox{if}\ i<n,\\
r_{i,i}^2=(2r_{i,i}+1)r_{j,j}, & \mbox{if}\ (i,j)\in\mathcal I_A,\\
(2r_{i,i}+1)\left(r_{i+1,n}+\sum\limits_{(i,j)\in\mathcal I_A}a_{i,j}r_{j,n}\right)=0, & \mbox{if}\ i<n-2,\\
(2r_{n-2,n-2}+1)r_{n-1,n}=0.
\end{array}
\right.
\end{equation}
Let us define the following function
\begin{equation}\label{func}
f(x)=\frac{x^2}{2x+1},\ \ \ \ x\in\mathbb K.
\end{equation}
By $Per_m(f)$ we denote a set of all $m$ periodic point of $f$, i.e.
$$
Per_m(f)=\left\{x\in\mathbb K: f^m(x)=x\right\},
$$
where $f^m=\underbrace{f\circ\dots\circ f}_{m}$

\begin{thm}\label{thmRB2}
Let ${\bf E}$ be an evolution algebra with structural matrix
$$
A=\left(
\begin{array}{llllllll}
0 & 1 & a_{1,3} & a_{1,4} & \dots & a_{1,n-1} & 0\\
0 & 0 & 1 & a_{2,4} &\dots & a_{2,n-1} & 0\\
0 & 0 & 0 & 1 & \dots & a_{3,n-1} & 0\\
\vdots & \vdots & \vdots & \vdots & \ddots & \vdots & \vdots\\
0 & 0 & 0 & 0 & \dots & 1 & 0\\
0 & 0 & 0 & 0 & \dots & 0 & 1\\
0 & 0 & 0 & 0 & \dots & 0 & 0
\end{array}\right)
$$ in a
natural basis
$\{{\bf e}_i\}_{i=1}^n$ and $R$ be a linear operator on $\bf E$. Then  $R\in\mathcal{RB}_1^\vartriangle({\bf E})$ if and only if
it has one of the following forms:
\begin{enumerate}
\item[$(i)$] if $\mathcal I_A=\emptyset$ then
$$
R({\bf e}_i)=\left\{
\begin{array}{ll}
\alpha{\bf e}_1+\beta{\bf e}_n, & i=1,\\
f^{i-1}(\alpha){\bf e}_i, & i>1
\end{array}
\right.\ \ \ \ \ \mbox{for some}\ \alpha\in\mathbb K\setminus\bigcup\limits_{k=0}^{n-2}f^{-k}(-\frac{1}{2}),\ \beta\in\mathbb K,
$$
here $f$ defined as \eqref{func}.
\item[$(ii)$] if $\mathcal I_A\neq\emptyset$ and
$i_0:=\left\{i: \exists j\ \mbox{such that}\ (i,j)\in\mathcal I_A\right\}$, then
\begin{equation}\label{RB119}
R({\bf e}_i)=\left\{
\begin{array}{ll}
f^{-i_0}(\alpha){\bf e}_1+\beta{\bf e}_n, & i=1,\\
f^{i-i_0-1}(\alpha){\bf e}_i, & i>1
\end{array}
\right.\ \ \ \ \ \mbox{for some}\ \alpha\in Per_m(f),\ \beta\in\mathbb K,
\end{equation}
here $I_1:=\left\{j-i-1: (i,j)\in\mathcal I_A\right\}$,
$$
I_2=\left\{|k-j|: k\neq j, \exists i\ \mbox{such that}\ (i,j),(i,k)\in\mathcal I_A\right\},
$$
and $m=GCD(I_1\cup I_2)$.
\end{enumerate}
\end{thm}

\begin{proof} $(i)$ Let $\mathcal I_A=\emptyset$. Then \eqref{RB1110} can be written as folllows
\begin{equation}\label{RB11100}
\left\{
\begin{array}{ll}
r_{i,i}^2=(2r_{i,i}+1)r_{i+1,i+1}, & \mbox{if}\ i<n,\\
(2r_{i,i}+1)r_{i+1,n}=0, & \mbox{if}\ i<n-1.
\end{array}
\right.
\end{equation}
Then the first equalities of \eqref{RB11100}
imply either $\prod_{i=1}^nr_{i,i}\neq0$ or
$$
r_{1,1}=r_{2,2}=\dots=r_{n,n}=0.
$$
If $r_{1,1}=r_{2,2}=\dots=r_{n,n}=0$ then from the second equalities of \eqref{RB11100}
we find that $r_{i,n}=0$ for all $i>1$. Hence, the Rota-Baxter operator
has the following form
\begin{equation}\label{R1}
R({\bf e}_i)=\left\{
\begin{array}{ll}
\beta{\bf e}_n, & i=1,\\
{\bf 0}, & i>1,
\end{array}
\right.\ \ \ \ \ \mbox{for some}\ \beta\in\mathbb K.
\end{equation}
Assume that $\prod_{i=1}^nr_{i,i}\neq0$. It is obvious that $r_{i,i}\neq-\frac{1}{2}$ for all $i<n$. Then, from
\eqref{RB11100} we get $r_{i,n}=0$ for every $i\in\{2,\dots,n-1\}$. Hence, the Rota-Baxter operator can be written as
\begin{equation}\label{R2}
R({\bf e}_i)=\left\{
\begin{array}{ll}
\alpha_1{\bf e}_1+\beta{\bf e}_n, & i=1,\\
\alpha_i{\bf e}_i, & i>1
\end{array}
\right.\ \ \ \ \ \mbox{for some}\ \alpha_1,\beta\in\mathbb K,
\end{equation}
here $\alpha_i=\frac{\alpha_{i-1}^2}{2\alpha_{i-1}+1}$, $i=\overline{2,n}$ such that $\alpha_i\notin\{-\frac{1}{2},0\}$.
By \eqref{R1} and \eqref{R2} we conclude that the Rota-Baxter operator has the form \eqref{RB119}.

$(ii)$ Let $\mathcal I_A\neq\emptyset$.
By the first equalities of \eqref{RB1110} we infer that $r_{i,i}\neq-\frac{1}{2}$. Then from
the last equality of \eqref{RB1110} we get $r_{n-1,n}=0$.
Putting it into the third equality of \eqref{RB1110} for $i=n-3$ we find $r_{n-2,n}=0$.
Assume that $r_{i_0+1,n}=r_{i_0+2,2}=\dots=r_{n-1,n}=0$ for some $2<i_0<n-1$.
Then putting these into the third equality of \eqref{RB1110} for $i=i_0-1$ and keeping in mind
$r_{i_0-1,i_0-1}\neq-\frac{1}{2}$ we obtain $r_{i_0,n}=0$. Thus, 
\begin{equation}\label{per}
R({\bf e}_i)=
\left\{
\begin{array}{ll}
\alpha_1{\bf e}_1+\beta{\bf e}_n, & \mbox{if}\ i=1,\\
\alpha_i{\bf e}_i, & \mbox{if}\ i>1,
\end{array}
\right.
\end{equation}
where $\alpha_1,\beta\in\mathbb K$ and
\begin{equation}\label{RB111000}
\left\{
\begin{array}{ll}
\alpha_i^2=(2\alpha_i+1)\alpha_{i+1}, & \mbox{if}\ i<n,\\
\alpha_i^2=(2\alpha_i+1)\alpha_{j}, & \mbox{if}\ (i,j)\in\mathcal I_A.
\end{array}
\right.
\end{equation}
It is easy to see that $\alpha_1=\dots=\alpha_n=0$ or $\alpha_1=\dots=\alpha_n=-1$ satisfies \eqref{RB111000}.
Assume that there exists $i_0<n-2$ such that $(i_0,i_0+2)\in\mathcal I_A$.
Then from \eqref{RB111000} we have
$\alpha_{i_0+1}=\alpha_{i_0+2}$. Consequently, after rewriting \eqref{RB111000} for $i=i_0+1$
one finds
$$
\alpha_{i_0+1}^2=(2\alpha_{i_0+1}+1)\alpha_{i_0+1},
$$
which holds only for $\alpha_{i_0+1}=0$ or $\alpha_{i_0+1}=-1$. Hence, by the first equalities of \eqref{RB111000}
we find $\alpha_{1}=\dots=\alpha_{n}=0$ or $\alpha_{1}=\dots=\alpha_{n}=-1$.
Thus, we have shown that $\alpha_{1}\in\{-1,0\}$ if $1\in I_1:=\left\{j-i-1: (i,j)\in\mathcal I_A\right\}$.
Let $\alpha_{1}=\alpha$ and $\alpha_{i}=f^{i-1}(\alpha)$, $i=\overline{2,n}$.
Assume that $(p,p+k+1)\in\mathcal I_A$. Due to $\alpha_{p}\neq-\frac{1}{2}$,
from \eqref{RB111000} we find $\alpha_{p+1}=\alpha_{p+k+1}$. On the other hand, we have
$\alpha_{p+k}=f^{k}(\alpha_{p})$. This means that $\alpha_{p+1}$ is $k$ periodic point of the function
$f$. Moreover, one has $\alpha_{i}=f^{i-p-1}(\alpha_{p+1})$, $i=\overline{p+1,n}$.
Let $(p',p'+k'+1),(p'',p''+k''+1)\in\mathcal I_A$. For the sake of convenience, assume that $p'\leq p''$.
Then
$$
\alpha_{i}=f^{i-p'-1}(\alpha_{p+1}),\ \ i=\overline{p'+1,n},
$$
where $\alpha_{p'+1}\in Per_{k'}(f)$.
This implies that
$$
\{\alpha_{p''+1},\dots,\alpha_n\}\subset\left\{\alpha_{p'+1},f(\alpha_{p'+1}),\dots,f^{k'-1}(\alpha_{p'+1})\right\}.
$$
On the other hand, we have
$$
\{\alpha_{p''+1},\dots,\alpha_n\}\subset\left\{\alpha_{p''+1},f(\alpha_{p''+1}),\dots,f^{k''-1}(\alpha_{p''+1})\right\}.
$$
Consequently, we infer that $\alpha_{p'+1}$ is a $k$ periodic point of $f$, where $k=GCD(k',k'')$.
We define
$$
I_2=\left\{|k-j|: k\neq j, \exists i\ \mbox{such that}\ (i,j),(i,k)\in\mathcal I_A\right\}.
$$
Note that $I_2$ can be empty. Assume that $I_2\neq\emptyset$, i.e. there exist
$(p,q),(p,r)\in\mathcal I_A$ such that $q<r$. Then from \eqref{RB111000} we find $\alpha_{p+1}=\alpha_{q}=\alpha_{r}$.
This means that $\alpha_{p+1}$ is $m$ periodic point of $f$, where $m=GCD(q-p,r-p)$.
Since $\mathcal I_A\neq\emptyset$ we can define the numbers
$$
i_0:=\left\{i: \exists j\ \mbox{such that}\ (i,j)\in\mathcal I_A\right\},\ \ \
m:=GCD(I_1\cup I_2).
$$
Then we conclude that \eqref{per} is Rota-Baxter operator if and only if $\alpha_{i_0+1}$ is $p$ periodic point
of $f$ and
$\alpha_i=f^{i-i_0-1}(\alpha_{i_0+1})$ for any $i\leq n$. This completes the proof.
\end{proof}
Let $\mathbb K=\mathbb R$. In this case the function \eqref{func} is regular, i.e.
for any $x\neq-\frac{1}{2}$ one has
$$
\lim\limits_{m\to\infty}f^m(x)=\left\{
\begin{array}{rr}
-1, & \mbox{if}\ \ i<-\frac{1}{2},\\[2mm]
0, & \mbox{if}\ \ i>-\frac{1}{2}.
\end{array}
\right.
$$
Then for the real case, Theorem \ref{thmRB2} can be formulated as follows

\begin{thm}\label{thmRB3}
Let ${\bf E}$ be a real evolution algebra with structural matrix
$$
A=\left(
\begin{array}{llllllll}
0 & 1 & a_{1,3} & a_{1,4} & \dots & a_{1,n-1} & 0\\
0 & 0 & 1 & a_{2,4} &\dots & a_{2,n-1} & 0\\
0 & 0 & 0 & 1 & \dots & a_{3,n-1} & 0\\
\vdots & \vdots & \vdots & \vdots & \ddots & \vdots & \vdots\\
0 & 0 & 0 & 0 & \dots & 1 & 0\\
0 & 0 & 0 & 0 & \dots & 0 & 1\\
0 & 0 & 0 & 0 & \dots & 0 & 0
\end{array}\right)
$$ in a
natural basis
$\{{\bf e}_i\}_{i=1}^n$. Then the following statements hold:
\begin{enumerate}
\item[$(i)$] if $\mathcal I_A=\emptyset$ then
$$
\mathcal{RB}_1^\vartriangle({\bf E})=
\left\{\left(
\begin{array}{lllllllllllll}
\alpha & 0 & 0 \dots & 0 & \beta\\
0 & f(\alpha) & 0 \dots & 0 & 0\\
\vdots & \vdots & \ddots & \vdots & \vdots\\
0 & 0 & 0 \dots & f^{n-2}(\alpha) & 0\\
0 & 0 & 0 \dots & 0 & f^{n-1}(\alpha)
\end{array}
\right): \alpha,\beta\in\mathbb R,\ \alpha\neq-\frac{1}{2}\right\}.
$$
\item[$(ii)$] if $\mathcal I_A\neq\emptyset$ then
$$
\mathcal{RB}_1^\vartriangle({\bf E})=
\left\{\left(
\begin{array}{lllllllllllll}
\alpha & 0 & 0 \dots & 0 & \beta\\
0 & \alpha & 0 \dots & 0 & 0\\
\vdots & \vdots & \ddots & \vdots & \vdots\\
0 & 0 & 0 \dots & \alpha & 0\\
0 & 0 & 0 \dots & 0 & \alpha
\end{array}
\right): \alpha=\pm1, \beta\in\mathbb R\right\}.
$$
\end{enumerate}
\end{thm}
\subsection{Case $\mathcal{RB}_1^\blacktriangle({\bf E})$}
The conditions \eqref{RB111} and \eqref{RB211}
are equivalent to the following
\begin{equation}\label{RB1111}
\left\{
\begin{array}{ll}
r_{i,i}^2=(2r_{i,i}+1)r_{i+1,i+1}, & \mbox{if}\ i<n,\\
r_{i,n-1}^2=(2r_{i,i}+1)\left(r_{i+1,n}+\sum\limits_{(i,j)\in\mathcal I_A}a_{i,j}r_{j,n}\right), & \mbox{if}\ i<n-1,\\
r_{i,m-1}^2+\sum\limits_{k\geq i:\atop{(k,m)\in\mathcal I_A}}a_{k,m}r_{i,k}^2=(2r_{i,i}+1)\left(r_{i+1,m}+\sum\limits_{(i,k)\in\mathcal I_A}a_{i,k}r_{k,m}\right), & \mbox{if}\ i+1<m<n,\\
r_{i,m-1}r_{j,m-1}+\sum\limits_{k>j:\atop{(k,m)\in\mathcal I_A}}a_{k,m}r_{i,k}r_{j,k}=r_{i,j}\left(r_{j+1,m}+\sum\limits_{(j,k)\in\mathcal I_A}a_{j,k}r_{k,m}\right), & \mbox{if}\ i<j<m<n,\\
r_{i,n-1}r_{j,n-1}=r_{i,j}\left(r_{j+1,n}+\sum\limits_{(j,k)\in\mathcal I_A}a_{j,k}r_{k,n}\right), & \mbox{if}\ i<j<n-2,\\
r_{i,n-2}r_{n-2,n-2}=r_{i,n-2}r_{n-1,n-1}, & \mbox{if}\ i<n-2,\\
r_{i,n-1}r_{n-2,n-1}=r_{i,n-2}r_{n-1,n}, & \mbox{if}\ i<n-2,\\
r_{i,n-1}r_{n-1,n-1}=r_{i,n-1}r_{n,n}, & \mbox{if}\ i<n-1.
\end{array}
\right.
\end{equation}
Let us rewrite the fourth equalities of \eqref{RB1111}
for $m=j+1$. Due to $a_{k,j+1}=0$ for every $k>j$, and $r_{k,j+1}=0$ for any $k>j+1$
we obtain
$$
r_{i,j}r_{j,j}=r_{i,j}r_{j+1,j+1},\ \ \ \ i<j<n-1.
$$
The last one together with the eighth equalities of \eqref{RB1111} yields that
\begin{equation}\label{ni1}
r_{i,j}r_{j,j}=r_{i,j}r_{j+1,j+1},\ \ \ \ i<j<n.
\end{equation}
By assumption (i.e. $R\in\mathcal{RB}_1^\blacktriangle({\bf E})$) one can find a pair $(i_0,j_0)$
such that $1\leq i_0<j_0<n$ and $r_{i_0,j_0}\neq0$. Then from \eqref{ni1}
we have $r_{j_0,j_0}=r_{j_0+1,j_0+1}$. Hence, by the first equalities of \eqref{RB1111}
we obtain either
$r_{1,1}=\dots=r_{n,n}=-1$ or $r_{1,1}=\dots=r_{n,n}=0$.

{\bf Case $\mathcal I_A=\emptyset$.}
Let us assume that $r_{1,i}=\beta_i$, $i>1$. Then from the second and the third
equalities of \eqref{RB1111} we find that
$$
r_{i,j}=
\left\{
\begin{array}{ll}
\beta_{j-i+1}^{2^{i-1}}, & \mbox{if}\ r_{1,1}=0,\\
-\beta_{j-i+1}^{2^{i-1}}, & \mbox{if}\ r_{1,1}=-1,
\end{array}
\right.\ \ \ \ \ \ 1<i<j.
$$
The fourth, the fifth and the seventh equalities of \eqref{RB1111} can be written as follows
\begin{equation}\label{457}
\left\{
\begin{array}{ll}
\beta_{m-i}^{2^{i-1}}\beta_{m-j}^{2^{j-1}}=\beta_{j-i+1}^{2^{i-1}}\beta_{m-j}^{2^{j}}, & \mbox{if}\ i<j<m-1<n-1,\\[2mm]
\beta_{n-i}^{2^{i-1}}\beta_{n-j}^{2^{j-1}}=\beta_{j-i+1}^{2^{i-1}}\beta_{n-j}^{2^{j}}, & \mbox{if}\ i<j<n-2,\\[2mm]
\beta_{n-i}^{2^{i-1}}\beta_{2}^{2^{n-3}}=\beta_{n-i-1}^{2^{i-1}}\beta_{2}^{2^{n-2}}, & \mbox{if}\ i<n-2.
\end{array}
\right.
\end{equation}
If $\beta_2\neq0$
$$
|\beta_m|=\beta_2^{2^{m-1}-1},\ \ 2<m<n.
$$
If $\beta_2=\dots=\beta_k=0\neq\beta_{k+1}$, $2\leq k<n-1$
$$
|\beta_m|=\left\{
\begin{array}{ll}
0, & m\neq1(\operatorname{mod}k),\\[2mm]
\beta_{k+1}^{\frac{2^{m-1}-1}{2^k-1}}, & m\equiv1(\operatorname{mod}k),
\end{array}
\right.\ \ \ \ k+1<m<n.
$$
\
{\bf Case $\mathcal I_A\neq\emptyset$.} In this setting, we have
$$
r_{i,i+1}=\left\{
\begin{array}{ll}
r_{1,2}^{2^{i-1}}, & \mbox{if}\ r_{1,1}=0,\\
-r_{1,2}^{2^{i-1}}, & \mbox{if}\ r_{1,1}=-1.
\end{array}
\right.\ \ \ \ 1<i<n.
$$

\section*{Acknowledgments}
	The authors thank the UAEU UPAR Grant No. G00004962 for support.

\section*{Declaration} The author declare that they have no conflict of~interests.

\section*{Availability of data and material} Not applicable.

\end{document}